\documentclass[paper]{paper_cls}
\usepackage[pdftex]{graphicx,xcolor}
\usepackage[fleqn]{amsmath}
\usepackage{newtxtext}
\usepackage[varg]{newtxmath}
\usepackage{bm}
\usepackage{comment}
\usepackage{cite}
\newtheorem{Theo}{Theorem}
\newtheorem{Assu}{Assumption}
\usepackage{booktabs}

\field{A}
\vol{109}
\no{xxx}
\title[1D Nonnegative Spline Smoothing via CSIP with a CP Method]{One-Dimensional Nonnegative Spline Smoothing via\\Convex Semi-Infinite Programming with a Cutting-Plane Method}
\authorlist{\authorentry[hiroki2580arai@keio.jp]{Hiroki ARAI}{m}{labelA}\MembershipNumber{}
 \authorentry[d.kitahara@keio.jp (Corresponding author)]{Daichi KITAHARA}{m}{labelA}\MembershipNumber{1616838}}
\affiliate[labelA]{The authors are with the
\EICdepartment{Dept.~of Applied Physics and Physico-Informatics}
\EICorganization{Keio University}
\EICaddress{Yokohama-shi, 223-8522 Japan}
}

\received{2026}{4}{xx}
\revised{xxxx}{4}{xx}
\onlinefirst{xxxx}{4}{xx}

\begin{document}
\maketitle
\begin{summary}
Spline functions are smooth piecewise polynomials widely used for interpolation and smoothing, and nonnegative spline smoothing is also studied for nonnegative data.
Previous research used sufficient conditions for the nonnegativity of spline functions because necessary and sufficient conditions for the nonnegativity are infinitely many linear inequalities, which are difficult to handle in optimization algorithms.
This conventional method quickly computes a nonnegative spline function via quadratic programming (QP), but the optimal solution may be slightly degraded by using the sufficient condition instead of the necessary and sufficient condition.
In this paper, we express one-dimensional (1D) nonnegative spline smoothing as a convex semi-infinite programming (CSIP) problem that directly deals with infinite inequality constraints.
As optimization algorithms for general SIP problems including nonconvex cases, local-reduction-based sequential quadratic programming (LRSQP) methods are used, but their~convergence performance deteriorates for certain problems due to multiple approximations during updates.
To quickly obtain the solution, we propose~a~cutting-plane (CP) method for the 1D nonnegative spline smoothing formulated as the CSIP problem.
In the proposed method, after giving an initial solution by the standard spline smoothing, we find the minimizer of each polynomial piece by using the closed-form solution for a low-degree polynomial or a numerical solution for a high-degree polynomial.
If the minimum value is negative, then such minimizer is added into the constraint of the problem to guarantee the nonnegativity.
This constrained problem is quickly solved~via QP, and we find the minimizer of each polynomial piece again.
We repeat these procedures until there are no negative minimum values.
The proposed method guarantees convergence to the original CSIP solution with (almost) no approximate computations, and its effectiveness is demonstrated in numerical experiments by comparison to the conventional methods, QP under the sufficient condition and CSIP using the MATLAB LRSQP algorithm.
\end{summary}


\section{Introduction}
Spline functions are continuous piecewise polynomials that are differentiable up to a finite order \cite{schumaker2007spline,Wahba-90,Ramsay-Silverman-05,reinsch1967smoothing,Chui-88}.
Due to their~ease of handling and certain optimality for the smoothness, spline functions are widely used in many fields.
In particular, they play an important role for data interpolation and smoothing, where interpolation refers to constructing a smooth function that passes through all noise-free data points, while smoothing refers to estimation of a noise-free smooth function from its noisy samples.
For both one-dimensional (1D) interpolation and smoothing, spline functions possess the optimality of minimizing the curvature, evaluated by the energy of the second derivative, among all $C^2$ functions, and the optimal functions are \textit{natural cubic splines} \cite{Wahba-90,Ramsay-Silverman-05,reinsch1967smoothing}.
Such a smoothing result is called \textit{1D spline smoothing} in this paper.

In mathematics and engineering, smooth functions that always take nonnegative values, such as probability density functions \cite{parzen1962estimation} and power spectral densities \cite{percival1993spectral}, are often~required.
When estimating such a nonnegative function~from nonnegative data points, the standard spline smoothing does not guarantee the nonnegativity over the entire domain.
In order to resolve this issue, \textit{nonnegative spline smoothing} is studied \cite{kitahara2015probability,kitahara2016two,kitahara2016thesis}, where sufficient conditions to guarantee the nonnegativity of spline functions are derived and used.
This nonnegative spline smoothing is formulated as a \textit{quadratic programming (QP) problem} with respect to the coefficient vector of a spline function, which is quickly and accurately solved by using an interior-point method \cite{Vanderbei-Carpenter-93,Altman-Gondzio-99}.
However, since the sufficient condition is imposed instead of a necessary and sufficient condition, the feasible set is smaller than the original feasible set of all nonnegative spline functions, which may result in a slightly degraded solution whose cost function value is higher than the original optimal solution.

In this paper, we aim to compute the exact solution~of 1D nonnegative spline smoothing without using the sufficient conditions in \cite{kitahara2015probability,kitahara2016two,kitahara2016thesis}.
We directly impose the constraint~that each polynomial piece of a spline function is nonnegative at every point in its corresponding interval.
This necessary~and sufficient condition for the nonnegativity is expressed as an infinite number   of linear inequalities, which are difficult to handle in standard optimization algorithms.
An optimization problem involving infinitely many inequality constraints is called a \textit{semi-infinite programming (SIP) problem} \cite{hettich1993semi,reemtsen1998semi,Goberna-Lopez-01,Lopez-Still-07,OkunoFukushima14,reemtsen1991discretization,blankenship1976infinitely,oustry2025convex,schwientek2021transformation}.
In particular, it is called a \textit{convex SIP (CSIP) problem} if~the cost function and the feasible set are convex.
Both in the 1D nonnegative spline smoothing are convex, and we can~construct a CSIP problem with respect to the coefficient vector.

As a SIP algorithm applicable even in nonconvex cases, a \textit{local-reduction-based sequential quadratic programming (LRSQP) method} \cite{reemtsen1998semi,OkunoFukushima14} is used.
In the LRSQP method,~the original infinite inequality constraints are replaced with finite inequality constraints in the vicinity of the current solution, and such locally reduced problems are approximately solved in a manner similar to the standard SQP method.
The LRSQP method is highly versatile, but its convergence performance deteriorates for certain problems since update~formulas~use some approximations that may decrease speed and accuracy.

In order to obtain the exact solution of the nonnegative spline smoothing more efficiently than the LRSQP method, we propose a \textit{cutting-plane (CP) method}.
CP methods were originally developed to solve optimization problems involving finite but complicated constraints \cite{integer,comp_nonconvex}.
Recently, for certain CSIP problems, a CP method that guarantees convergence to the optimal solution was proposed \cite{oustry2025convex}.
Our~proposed method is based on this CP method and specifically~designed for the 1D nonnegative spline smoothing as CSIP.

After obtaining an initial solution by the standard spline smoothing, the proposed method finds the minimizer of each polynomial piece by using the closed-form solutions for low-degree cases and numerical solutions for high-degree cases.
If there are negative minimum values, those minimizers are added into the constraint to guarantee the nonnegativity.
This constrained problem is QP and quickly solved by an interior-point method.
These procedures are repeated until there are no negative minimum values.
The proposed update formulas have (almost) no approximation errors, differently from the LRSQP method.
In numerical experiments, we show that~the proposed CP method and the LRSQP method provide more accurate solutions than the conventional QP method with the sufficient condition, and the proposed method achieves~much shorter computation times and slightly better accuracy than the LRSQP method.
The superiority of the closed-form~solutions to numerical solutions in the minimizer search of~the proposed method is also demonstrated by the experiments.

This paper is organized as follows.
Sections 2.1 and 2.2 introduce the definition of spline functions and the standard spline smoothing.
Section 2.3 explains the conventional~nonnegative spline smoothing with the sufficient condition \cite{kitahara2016thesis}. 
Section 3.1 formulates the nonnegative spline smoothing as a CSIP problem, and Section 3.2 explains the LRSQP method usable in MATLAB\@. 
Section 4.1 summarizes the CP method in \cite{oustry2025convex}, and Section 4.2 presents detailed calculations of the proposed method.
Sections 5.1 and 5.2 explain results~of~the numerical experiments, and Section 6 concludes this paper.

\section{Mathematical Preliminaries}
Let $\mathbb R$ and $\mathbb N$ be the sets of all real numbers and nonnegative integers, respectively.
For any nonnegative integer $\rho\in\mathbb N$~and any two real numbers $\alpha,\beta\in\mathbb R$ s.t.~$\alpha<\beta$, $C^{\rho}(\alpha,\beta)$ denotes the set of all $\rho$-times continuously differentiable functions~on $[\alpha,\beta]$.
For any $d\in\mathbb N$, $\mathbb P_{d}$ stands for the set of all univariate polynomials of degree $d$ (at most), i.e., $\mathbb P_{d}:=\{p:\mathbb R\ni x\mapsto \sum^{d}_{k=0} c_{k}x^{k}\in\mathbb R \,|\,c_k \in\mathbb R\}$.
We write vectors and matrices by boldface small and capital letters, respectively.

\subsection{Definition and Representation of Spline Functions}
Let ${\sqcup}_m:=\{[\xi_{i-1},\xi_i]\}_{i=1}^{m}$ be a set of $m$ subintervals $[\xi_{i-1},\xi_i]$ divided by \textit{knots} $\xi_{i}$ ($\xi_0 < \xi_1 < \cdots < \xi_m$).
For~a partition~${\sqcup}_m$ and two integers $\rho, d \in \mathbb N$ s.t.~$0 \le \rho < d$, define the set of all spline functions of degree $d$ and smoothness $\rho$ on ${\sqcup}_m$ by
\mathindent=2.7mm
\[
\mathcal S^{\rho}_{d}(\sqcup_m):=\{s\in C^{\rho}(\xi_0,\xi_m)\,|\,s=p_{i}\in\mathbb P_{d}\mbox{ on }[\xi_{i-1},\xi_i]\}\mbox{.}
\]
In particular, for $k \in \mathbb N$, a spline function $s \in \mathcal S^{2k}_{2k+1}(\sqcup_m)$ is called a \emph{natural spline} of degree $2k+1$ if its $l$th derivatives $s^{(l)}$ ($l=k+1,k+2 \ldots, 2k$) satisfy $s^{(l)}(\xi_0)=s^{(l)}(\xi_m)=0$.

In this paper, for $s\in \mathcal S^{\rho}_{d}(\sqcup_m)$, its $i$th polynomial piece $p_i$ on $[\xi_{i-1},\xi_i]$ is expresses as the \textit{Bernstein-B{\'e}zier form}
\mathindent=0.9mm
\begin{align}
p_i(x)&= \sum_{k=0}^{d} b_{d(i-1)+k}\,\frac{d!}{(d-k)!k!}\,
\Bigl(\frac{\xi_i - x}{\Delta_i}\Bigr)^{d-k}
\Bigl( \frac{x - \xi_{i-1}}{\Delta_i} \Bigr)^k \nonumber\\
&= \sum_{k=0}^{d} b_{d(i-1)+k}\,
\frac{d!}{(d-k)!k!}\,
(1-\tau)^{d-k} \tau^k\mbox{,}\label{Bernstein}
\end{align}
where $b_{d(i-1)+k} \in \mathbb R$ ($k=0,1,\ldots,d$) are the coefficients of each polynomial $p_i \in \mathbb P_d$ ($i=1,2,\ldots,m$), $\Delta_i := \xi_i - \xi_{i-1}$, and $\tau := (x-\xi_{i-1})/\Delta_i$.
In each subinterval, the value of $\tau$~is normalized, i.e., $\tau\in[0,1]\Leftrightarrow x\in[\xi_{i-1}, \xi_i]$.
In addition,~the function value at each knot $\xi_i$ is easily given by $s(\xi_i)=b_{di}$.

\subsection{Standard Spline Smoothing and Its QP Formulation}
We consider the estimation of a continuous function $y$ from its finite samples $y_i = y(x_i) + \varepsilon_i$ including noise $\varepsilon_i\in\mathbb R$~at $x_0 < x_1 < \cdots < x_n$.
To accurately reconstruct the original smooth function $y$, it is necessary to find a function that is as smooth as possible while allowing a certain degree of error at each data point $(x_i, y_i)$.
This process is called smoothing, and in the 1D case, the following theorem is known \cite{reinsch1967smoothing,Wahba-90,Ramsay-Silverman-05}.

\begin{Theo}[Spline Smoothing]
Suppose $n \ge 1$ and $n+1$ data points $\{(x_i, y_i)\}_{i=0}^{n}$ s.t.~$x_0 < x_1 < \cdots < x_n $ are given.
The following optimization problem (variational problem)
\mathindent=6mm
\begin{equation}
\mathop{\mathrm{minimize}}_{u \in C^{2}(x_0, x_n)}\, \sum_{i=0}^{n} |y_i - u(x_i)|^{2}+ \lambda \int_{x_0}^{x_n} |u''(x)|^{2} \, \mathrm{d}x \label{smoothing}
\end{equation}
has the \textit{unique} optimal solution for any $\lambda>0$.
It is a~\textit{natural cubic spline} $s \in \mathcal{S}^{2}_{3}(\sqcup_{n})$ with knots $\xi_i := x_i$ ($i = 0, 1, \ldots, n$).
\end{Theo}

The first term in \eqref{smoothing} represents the squared error at each data point, while the second term evaluates the roughness penalty for $u$ by the energy of the second derivative $u''$.
The weight $\lambda>0$ for the second term is called the \textit{smoothing~parameter} and controls the trade-off between the data fidelity and the smoothness.
Since the solution is a spline function, the problem in \eqref{smoothing} is called \textit{1D spline smoothing}.
By~computing the optimal spline function with an appropriate $\lambda$, the smooth function $y$ can be reconstructed with high accuracy.

When performing the spline smoothing in (\ref{smoothing}), the function space can be restricted from $C^{2}(x_0, x_n)$ to $\mathcal{S}^{2}_{3}(\sqcup_n)$, and the optimal natural cubic spline can be quickly computed for a fixed $\lambda$ by closed-form formulas.
To connect the discussion with nonnegative spline smoothing in the next section, in the following, we give a QP formulation applicable not only to $\mathcal{S}^{2}_{3}(\sqcup_n)$ but also to general spline spaces $\mathcal{S}^{2}_{d}(\sqcup_m)$ s.t.~$d\geq3$ and $[x_0, x_n]\subset[\xi_0, \xi_m]$, rather than the closed-form solution.

The coefficient vector $\bm b = (b_0, b_1, \ldots, b_{dm})^{\mathrm T}\in \mathbb{R}^{dm+1}$ is regarded as the parameter that determines the shape of a spline function $s \in \mathcal{S}^2_d(\sqcup_m)$.
From (1), for the first term~in (2), we can construct some sparse matrix $\bm A \in \mathbb{R}^{(n+1)\times (dm+1)}$ satisfying $\bm A \bm b = (s(x_0), s(x_1), \ldots, s(x_n))^\mathrm{T}\in \mathbb{R}^{n+1}$.
For the second term with the integral interval $[\xi_0, \xi_m]$, we have
\mathindent=1mm
\[
\begin{aligned}
&\int_{\xi_0}^{\xi_m}|s''(x)|^2\,\mathrm{d}x
=\sum_{i=1}^{m}\int_{\xi_{i-1}}^{\xi_i} |p_i''(x)|^2 \, \mathrm{d}x\\
&\quad= \sum_{i=1}^{m} \sum_{k=0}^{d} \sum_{l=0}^{d} b_{d(i-1)+k} b_{d(i-1)+l}\,\Biggl(\frac{1}{\Delta^3_{i}}\sum_{v=0}^{2} \sum_{w=0}^{2}
e_v e_w q_{v,w}^{k,l}\Biggl) \\
&\quad=\sum_{i=1}^m\bm b_{\langle i\rangle}^{\mathrm{T}}\bm Q_{\langle i\rangle} \bm b_{\langle i\rangle}=\bm b^\mathrm{T}\bm Q\bm b\mbox{,}
\end{aligned}
\]
where $q_{v,w}^{k,l}:=\frac{(d!)^2(2d-k-l+v+w-4)!(k+l-v-w)!}{(2d-3)!(d-k+v-2)!(k-v)!(d-l+w-2)!(l-w)!}$ if $k\in[v, d+v-2]$ and $l\in[w, d+w-2]$, and $q_{v,w}^{k,l}:=0$ otherwise, with $e_0=1$, $e_1:=-2$, and $e_2:=1$ (see \cite[Proposition 7]{kitahara2016thesis}).
$\bm b_{\langle i\rangle}=(b_{d(i-1)}, b_{d(i-1)+1},\ldots , b_{di})^{\mathrm{T}}\in \mathbb{R}^{d+1}$ is a subvector of $\bm{b}$ and the coefficient vector of the $i$th polynomial $p_i$, and $\bm Q_{\langle i\rangle}\in \mathbb{R}^{(d+1)\times{(d+1)}}$ and $\bm Q\in\mathbb{R}^{(dm+1)\times (dm+1)}$ are dense and sparse symmetric positive semidefinite matrices.
 
Since two adjacent polynomials $p_i$ and $p_{i+1}$ share one coefficient $b_{di}$ ($i=1, 2, \ldots, m-1$), $s \in C^0(\xi_0, \xi_m)$ always holds for every $\bm b\in\mathbb R^{dm+1}$. On the other hand, the condition $s \in C^2(\xi_0, \xi_m)$ is not guaranteed, and it is equivalent to the following linear equations (see \cite[Proposition 8]{kitahara2016thesis})
\mathindent=0.2mm
\[
\begin{aligned}
\frac{1}{\Delta_{i}^l}\sum_{k=0}^{l} \frac{(-1)^{l-k}l!}{(l-k)!k!}\,b_{di+k-l}
-\frac{1}{\Delta_{i+1}^l}\sum_{k=0}^{l} \frac{(-1)^{l-k}l!}{(l-k)!k!}\,b_{di+k}
=0\ \ \\
\mbox{(}l= 1, 2\; \mbox{and}\; i=1, 2, \ldots, m-1\mbox{).}
\end{aligned}
\]
We express these equations as a linear constraint $\bm H \bm b = \bm 0\Leftrightarrow s\in\mathcal S^{2}_{d}(\sqcup_m)$ with some sparse matrix $\bm H \in \mathbb{R}^{(2m-2)\times (dm+1)}$.

To summarize the discussion, the spline smoothing in (2) within $\mathcal S^{2}_{d}(\sqcup_m)$ is reduced to the following QP problem
\mathindent=3mm
\begin{equation}
\mathop{\mathrm{minimize}}_{\bm b \in \mathbb{R}^{dm+1}}
\ \lVert \bm y - \bm A \bm b \rVert_2^2
+ \lambda \bm b^{\mathrm T} \bm Q \bm b
\quad \mbox{subject to}\;  \bm H \bm b = \bm 0\mbox{,}\label{2d}
\end{equation}
where $\bm{y}:=(y_0,y_1,\ldots,y_n)^{\mathrm{T}}\in\mathbb R^{n+1}$ and $\lVert\cdot\rVert_{2}$ is the $\ell_{2}$-norm.
The cost function in (3) is \textit{strongly convex} on the constraint set (see Appendix A), and the \textit{unique} optimal solution can~be efficiently obtained by an interior-point algorithm \cite{Vanderbei-Carpenter-93,Altman-Gondzio-99}.

\subsection{Conventional Nonnegative Spline Smoothing}
Spline smoothing has also been studied for applications involving nonnegative data, such as the estimation of probability density functions and power spectral densities\cite{kitahara2015probability,kitahara2016two,kitahara2016thesis}.
Even if all data points $\{(x_i, y_i)\}_{i=0}^{n}$ satisfy $y_i \geq 0$, the overall nonnegativity of a natural cubic spline, which is the optimal solution of the problem in \eqref{smoothing}, is not guaranteed.
To ensure the overall nonnegativity of spline functions, it is necessary to explicitly introduce the nonnegativity constraint and solve the following optimization problem
\mathindent=7mm
\begin{align}
&\!\mathop{\mathrm{minimize}}_{s \in \mathcal{S}^2_d(\sqcup_m)}\,
\sum_{i=0}^{n} |y_i - s(x_i)|^2
+ \lambda \int_{\xi_0}^{\xi_m} |s''(x)|^2 \, \mathrm{d}x \nonumber\\
&\ \quad \mbox{subject to}\quad s(x) \ge 0\,\; \mbox{for all}\; x \in [\xi_0, \xi_m]\mbox{,}\label{NN}
\end{align}
which is referred to as \textit{1D nonnegative spline smoothing}.

In \eqref{NN}, the nonnegativity of $s$ has to be satisfied for all points $x$ in $[\xi_0, \xi_m]$, i.e., each polynomial piece $p_i$ has to be nonnegative over the entire interval $[\xi_{i-1}, \xi_i]$.
However, for $d\geq3$, it may be impossible to give a simple necessary~and sufficient condition that guarantees the nonnegativity of $p_i$ over $[\xi_{i-1}, \xi_i]$.
Therefore, previous research \cite{kitahara2015probability, kitahara2016two, kitahara2016thesis} derived a sufficient condition represented as a finite number of linear inequalities and used it instead of a necessary and sufficient condition.
For the Bernstein-B{\'e}zier  form in (1), the sufficient condition\footnote{This is equivalent to the sufficient condition in \cite[Eq.~(23)]{Yao-25} for the different spline representation $s(x)=p_{i}(x)=\sum_{k=0}^{d}c^{\langle i \rangle}_{k}\tau^{k}$.} for the nonnegativity of $p_i$ is simply given by
\mathindent=7mm
\begin{align}
&b_{d(i-1)+k} \geq 0\quad\mbox{($k = 0, 1, \ldots, d$)} \nonumber\\
&\qquad\Rightarrow\quad p_i(x) \geq 0 \,\; \mbox{for all}\; x \in [\xi_{i-1}, \xi_i].
\label{SC}
\end{align}
From (5), the sufficient condition for the nonnegativity of a spline function $s\in\mathcal S^{2}_{d}(\sqcup_m)$ is expressed as $\bm b\geq\bm 0$.
As~a~result, the nonnegative spline smoothing in (4) under the sufficient condition was formulated as the following QP problem
\mathindent=0.1mm
\begin{equation}
\mathop{\mathrm{minimize}}_{\bm b \in \mathbb{R}^{dm+1}}
\ \lVert \bm y - \bm A \bm b \rVert_2^2
+ \lambda\bm b^{\mathrm T} \bm Q \bm b
\quad \mbox{subject to}\; \bm b \ge \bm 0\mbox{, } \bm H \bm b = \bm 0\mbox{.}\label{SC2d}
\end{equation}
This is also efficiently solved by an interior-point method.

The (linear) inequality constraint in \eqref{SC2d} is not equivalent to the nonnegativity condition of the spline function $s$, but is a tractable sufficient condition.
Therefore, the set of spline functions corresponding to coefficient vectors in the feasible set $\{\bm{b}\in\mathbb R^{dm+1}\,|\, \bm b \ge \bm 0\mbox{ and }\bm H \bm b = \bm 0 \}$ is slightly smaller~than the ideal set of all nonnegative spline functions.
As a result, depending on the placement of the data points $\{(x_i, y_i)\}_{i=0}^{n}$, the 
optimal solution in \eqref{SC2d} may have a higher cost function value than the optimal nonnegative spline function in (4).

\section{Nonnegative Spline Smoothing via CSIP}
\subsection{CSIP Formulation}
In this paper, we find the optimal solution of the 1D nonnegative spline smoothing problem without relying on a sufficient condition such as (5) for the nonnegativity, i.e., we compute the solution of the problem in (4) without restricting the feasible set.
To this end, we directly impose the nonnegativity condition of each polynomial $p_i$ on its coefficient~vector~$\bm b_{\langle i\rangle}$.
This yields an infinite number of linear inequality constraints for the coefficient vector $\bm b$, and hence the problem becomes difficult to solve by standard optimization algorithms.

The class of optimization problems involving infinitely many inequality constraints is called \textit{semi-infinite programming (SIP)} \cite{hettich1993semi,reemtsen1998semi,Goberna-Lopez-01,Lopez-Still-07,OkunoFukushima14,reemtsen1991discretization,blankenship1976infinitely,oustry2025convex,schwientek2021transformation}, and the nonnegative spline smoothing in \eqref{NN} can be formulated as a \textit{convex SIP (CSIP)} problem.
This is because the cost function is a convex quadratic function on $\bm b$, and the feasible set is the intersection of convex sets, finite linear equality constraints and infinitely many linear inequality constraints.
In the following, we first describe the general form of SIP and then formulate the 1D nonnegative spline smoothing in \eqref{NN} as a CSIP problem.

In general, SIP is formulated as the following problem
\mathindent=4mm
\begin{align}
\mathop{\mathrm{minimize}}_{\bm{x}\in\mathbb R^{\zeta}}&
\,\, f(\bm{x}) \nonumber\\
\mbox{subject to}&\ \biggl\{\begin{aligned}
 &g_i(\bm x, \tau) \le 0\ \; \forall\tau \in T_i\quad\mbox{(}i = 1,2,\ldots,I\mbox{),} \\
& h_j(\bm x) = 0\quad\mbox{(}j = 1,2,\ldots,J\mbox{),}
\end{aligned}\label{SIP}
\end{align}
where $f:\mathbb{R}^\zeta\to\mathbb{R}$ is the cost function, $g_i:\mathbb{R}^{\zeta}\times T_i\to\mathbb{R}$ with $T_i:=[\alpha_{i}, \beta_{i}]$ s.t.~$\alpha_{i}<\beta_{i}$ ($i=1,2,\ldots,I$) are functions corresponding to infinitely many inequality constraints, and $h_j:\mathbb{R}^{\zeta}\to\mathbb{R}$ ($j=1,2,\ldots, J$) are functions corresponding to finite equality constraints.
The number of the functions~$g_i$ is finite, but each condition $g_i(\bm x,\tau)\le 0$ has to hold for all $\tau\in[\alpha_{i},\beta_{i}]$, which means that the number of the inequality constraints for the vector variable~$\bm x$ is infinite.
To solve~the SIP problem in (7), it is commonly assumed that $f$, $g_i$, and $h_j$ are differentiable up to the first or second order \cite{hettich1993semi,reemtsen1998semi}.

In this paper, for simplicity, we assume that every $T_i$~is normalized as $T_i=T=[0,1]$.
In particular, if $f$ is convex, each $g_i(\cdot,\tau)$ is also convex with respect to $\bm x$ for any $\tau\in T$, and each $h_j$ is affine, then the problem in \eqref{SIP} is called CSIP since both the cost function and the feasible set are convex.
Note that when $h_j$ is convex but not affine, the feasible set is generally nonconvex.
Thus, $h_j$ should be affine, not convex.

Next, we show that the nonnegative spline smoothing~in \eqref{NN} is regarded as a CSIP problem.
From (1), the nonnegativity condition of each polynomial $p_i$ can be expressed~as
\mathindent=2mm
\begin{align}
\sum_{k=0}^{d} \frac{d!(1-\tau)^{d-k} \tau^k}{(d-k)!k!}\,b_{d(i-1)+k} =\bm g_\tau^{\mathrm T}\bm b_{\langle i\rangle}  \geq 0
\quad\forall \tau \in T\mbox{,}\label{NNc}
\end{align}
where $\bm g_\tau:=\mathrm {vec}(\frac{d!(1-\tau)^{d-k} \tau^k}{(d-k)!k!})_{k=0}^{d}\in \mathbb R^{d+1}$ is defined from~$\tau$.
By adopting the coefficient vector $\bm{b}$ as the variable to be~optimized as with the problems in \eqref{2d} and \eqref{SC2d}, and by adding the constraint in \eqref{NNc} to the problem in \eqref{2d}, the nonnegative spline smoothing in \eqref{NN} is reduced to the following problem
\mathindent=4.4mm
\begin{align}
\mathop{\mathrm{minimize}}_{\bm{b}\in \mathbb{R}^{dm+1}}&
\,\, \lVert \bm{y}-\bm{A}\bm{b}\rVert_{2}^{2}
+ \lambda\bm{b}^{\mathrm{T}}\bm{Q}\bm{b} \nonumber\\
\mbox{subject to}&\ \biggl\{\begin{aligned}
 &\bm g_\tau^{\mathrm T}\bm b_{\langle i\rangle}  \geq 0\ \; \forall\tau \in T\quad\mbox{(}i = 1,2,\ldots,m\mbox{),} \\
& \bm h_j^{\mathrm T}\bm b = 0\quad\mbox{(}j = 1,2,\ldots,2m-2\mbox{),}
\end{aligned}\label{cSIP}
\end{align}
where $\bm h_j^{\mathrm T}$ is the $j$th row vector of the matrix $\bm H$ in (3) and (6).
By letting $\zeta:=dm+1$, $\bm x:=\bm b$, $I:=m$, $J:=2m-2$,~$f(\bm b):= \lVert \bm{y}-\bm{A}\bm{b}\rVert_{2}^{2}+\lambda\bm{b}^{\mathrm{T}}\bm{Q}\bm{b}$, $g_i(\bm b,\tau):={-}\bm g_\tau^{\mathrm T}\bm b_{\langle i\rangle}$ and $h_j(\bm b):=\bm h_j^{\mathrm T}\bm b$ in (7), we find that the problem in \eqref{cSIP} can be considered as a SIP problem.
Moreover, since the cost function $f$ is convex and the constraint functions $g_i(\cdot,\tau)$ and $h_j$ are linear on $\bm{b}$, the nonnegative spline smoothing in~(9) is a CSIP problem.

\subsection{MATLAB LRSQP Algorithm}
SIP includes infinitely many constraints, but it is very difficult to construct optimization algorithms that keep handling the infinite constraints in update formulas for general cases.
A key concept common to SIP algorithms is iteratively solving optimization problems of finite constraints that are given by the \textit{discretization-based relaxation} \cite{reemtsen1991discretization,schwientek2021transformation,OkunoFukushima14} as follows.

Let $P[\mathcal T]$ denote the original SIP problem in (7) with $\mathcal T:= \bigtimes^I_{i=1}T_i=[0,1]^I$.
For $P[\mathcal T]$, the SIP algorithms~generally consider the following finite-constraint relaxation
\mathindent=3mm
\begin{align}
\mathop{\mathrm{minimize}}_{\bm{x}\in\mathbb R^{\zeta}}&
\,\,  \tilde{f}(\bm{x}) \nonumber\\
\mbox{subject to}&\ \biggl\{\begin{aligned}
 & \tilde{g}_i(\bm x, \tau) \le 0\ \; \forall\tau \in \widetilde{T}_i\quad\mbox{(}i = 1,2,\ldots,I\mbox{),} \\
&  \tilde{h}_j(\bm x) = 0\quad\mbox{(}j = 1,2,\ldots,J\mbox{),}
\end{aligned}\label{rSIP}
\end{align}
and this problem is denoted by $P[\widetilde{\mathcal T}]$ with $\widetilde {\mathcal T}:=\bigtimes^I_{i=1}\widetilde T_i$.
In (10), each $\widetilde{T}_i\subset T$ is a \textit{discrete subset} of the~interval $[0,1]$.
For $P[\widetilde{\mathcal T}]$, since the number of the inequality constraints for the variable $\bm x$ is finite, we can find the optimal solution (or an approximate solution) by standard optimization algorithms.
Note that the functions $f$, $g_i$, and $h_j$ in (7) are also replaced with their approximations $\tilde f$, $\tilde g_i$, and $\tilde h_j$.
We can simply set $\tilde f=f$, $\tilde g_i=g_i$, and $\tilde h_j=h_j$, but several algorithms \cite{reemtsen1998semi,OkunoFukushima14} iteratively update approximate functions $\tilde f$, $\tilde g_i$, and~$\tilde h_j$ such~as $\tilde f_{r}$, $\tilde g_{i,r}$, and $\tilde h_{j,r}$, where $r\in\mathbb N$ denotes the update index.

As a baseline algorithm for comparison with the proposed method in Section 4, we introduce a \textit{local-reduction-based sequential quadratic programming (LRSQP)} method, implemented in MATLAB as \texttt{fseminf} function \cite{matlab_optim_doc}.
In~the LRSQP method, the cost function $f$ and the constraint functions $g_i$ and $h_j$ are approximated in a manner similar to the standard SQP method, where each $\widetilde{T}_i$ is \textit{locally given} in the neighborhood of the current solution as the set of all distinct local maximizers $\tau\in [0,1]$ of $g_i(\bm x,\cdot)$ (see \cite{reemtsen1998semi} for the definition of $\widetilde{T}_i$).
As a result, in the LRSQP method, $\tilde{f}$ becomes a quadratic function, and $\tilde g_i$ and $\tilde h_j$ become affine functions.

An~outline of the LRSQP method is described below.
\begin{itemize}
\item[\mbox{(i)}]
Fix $\bm x$ at an initial solution $\bm x_{-1}$ given as one input~argument.
Detect all $\tau \in T$ at which~$g_i(\bm x_{-1},\cdot)$~reaches local maxima, and define the set of the detected points~$\tau$ as $\widetilde{T}_{i,0}$ ($i=1,2,\ldots,I$).
Give an initial quadratic~approximation $\tilde f_{0}$ and initial affine approximations $\tilde g_{i, 0}$ and $\tilde h_{i, 0}$ around $\bm{x}_{-1}$.
Set $r=0$ and proceed to (ii).

\item[\mbox{(ii)}]
Solve the relaxed problem $P[\widetilde {\mathcal T}_{r}]$ with $\widetilde{\mathcal{T}}_{r}:=\bigtimes^I_{i=1}\widetilde T_{i,r}$, where $P[\widetilde {\mathcal T}_{r}]$ is \textit{QP with~respect to the direction vector}~$\bm d$ from $\bm x_{r-1}$, and its solution $\bm d_{r}$ is given by a QP solver.
Perform a line search along the direction $\bm{d}_{r}$ to obtain an approximate solution of the original problem $P[\mathcal T]$~as $\bm x_{r} = \bm x_{r-1} + a_{r} \bm d_r$ with a stepsize $a_{k}$.
Proceed to (iii).

\item[\mbox{(iii)}]
Fix $\bm x$ at the current solution $\bm x_r$ given in (ii).
Detect all $\tau \in T$ at which $g_i(\bm x_{r},\cdot)$~reaches local maxima.
Define the set of the detected $\tau$ as $\widetilde{T}_{i,r+1}$.
Proceed to (iv).

\item[\mbox{(iv)}]
Check the \textit{Karush--Kuhn--Tucker (KKT) conditions} for $P[\mathcal T]$.
If the KKT conditions hold within a tolerance, terminate the algorithm and return $\bm{x}_r$ as the solution of $P[\mathcal T]$.
Otherwise, increment $r$ by 1 and proceed to (v).

\item[\mbox{(v)}]
Update the quadratic function $\tilde f_r$ and the affine functions $\tilde g_{i, r}$ and $\tilde h_{i, r}$ around $\bm{x}_{r-1}$.
Proceed to (ii).
\end{itemize}

Differently from the standard SQP method, the LRSQP method (approximately) finds every local maximizer of each $g_i(\bm x, \cdot)$ with respect to $\tau\in T$ in Steps (i) and (iii).
Moreover, each local maximizer gradually changes during~the~iterative computations, i.e.,  the sampling points included in $\widetilde{T}_{i,r}$ differ for each $r$, and hence $\widetilde{T}_{i, r} \not\subset \widetilde{T}_{i,r+1}$ holds\footnote{If the algorithm converges, $\widetilde{T}_{i, r} $ and $\widetilde{T}_{i,r+1}$ become extremely close to each other when $r$ approaches infinity.} in general.
In the MATLAB LRSQP algorithm (\texttt{fseminf}), the~quadratic~function $\tilde f_{i,r}$ is constructed by the \textit{Broyden--Fletcher--Goldfarb--Shanno (BFGS) method} \cite{dai2013perfect} in Steps (i) and (v).
By default, gradient calculations in the BFGS method, the construction of $\tilde g_{i,r}$ and $\tilde h_{i,r}$, and the evaluation of the KKT conditions are performed by using \textit{numerical differentiation}.
In Step~(ii),~an \textit{active-set method} \cite{active} is used as a QP solver to compute the optimal solution $\bm{d}_{r}$ of $P[\widetilde{\mathcal T}_r]$.
In Step (i) and (iii), each local maximizer $\tau\in T$ is \textit{approximately} detected with the use of local cubic  (or quadratic) polynomial interpolation.

For the nonnegative spline smoothing in (9), $g_i(\bm b,\tau)={-}\bm g_\tau^{\mathrm T}\bm b_{\langle i\rangle}$ and $h_j(\bm b)=\bm h_j^{\mathrm T}\bm b$ are linear with respect to $\bm b$, and thus there are no approximation errors for $g_{i}$ and $h_{i}$ including their gradient calculations by numerical differentiation, i.e., $\tilde g_{i, r}=g_i$ and $\tilde h_{i, r}=h_i$ hold for all $r\in\mathbb N$.
On the~other hand, some approximation errors arise for the cost function $f$ and its gradient $\nabla f$ because of the BFGS method and numerical differentiation, even though $f(\bm b)= \lVert \bm{y}-\bm{A}\bm{b}\rVert_{2}^{2}+\lambda\bm{b}^{\mathrm{T}}\bm{Q}\bm{b}$ is a quadratic function from the beginning.
Furthermore,~since the local maximizers are detected through local cubic fittings in Steps (i) and (iii), errors also arise for the sampling~points in $\widetilde{T}_{i,r}$ for $d\geq 4$.
These errors are expected to decrease~the convergence speed and the numerical accuracy of the solution (see Appendix B for the KKT conditions in Step (iv)).

Some LRSQP implementations are guaranteed to converge from any~initial solution to a point satisfying the KKT conditions \cite{reemtsen1998semi,OkunoFukushima14}.  
However, it is very difficult to fully figure out the internal implementation of the MATLAB \texttt{fseminf} function, and it has not been clarified, to the best of the~authors' knowledge, which paper's implementation is the same. 
Therefore, in this paper, we experimentally verify its convergence speed and numerical accuracy in Section 5.
In the~next section, we propose another algorithm that is guaranteed to converge to the optimal solution of the CSIP problem in (9).

\section{CP Method for Fast and Accurate Optimization}
The MATLAB \texttt{fseminf} function is a general SIP solver for problems including nonconvex cost and nonlinear~constraint functions.
Although the application range of this algorithm~is wide, its update formulas use multiple approximations, such as replacing both cost and constraint functions and the local cubic fittings for detecting the local maximizers $\tau\in T$.
Thus, for certain SIP problems, there is still room for improvement in terms of computation time and numerical accuracy.

In this section, we propose an optimization algorithm, specifically designed for the nonnegative spline smoothing~in (9), that directly uses $f$, $g_i$, and $h_j$ without approximations and detects the \textit{global maximizer} of $g_i(\bm b, \cdot)$ with respect to $\tau\in T$.
The proposed algorithm and its convergence analysis are based on the \textit{cutting-plane (CP) method} in \cite{oustry2025convex}, and the optimal solution of the problem in (9) is computed more~efficiently and accurately compared to the LRSQP method.
In the following, we first describe the CP method in \cite{oustry2025convex} and its assumptions for the SIP problem in (7), and then present the proposed algorithm for the nonnegative spline smoothing.

\subsection{Convergence of a CP Method Containing Slight Errors}
CP methods were originally developed for optimization problems such as \textit{integer linear programming} \cite{integer} and \textit{nonlinear programming with complicated constraints} \cite{comp_nonconvex}.
They solve a relaxed problem after removing constraints that are difficult to handle explicitly.
If the solution of the relaxed problem~is not included in the original feasible set, solve a new relaxed problem with \textit{one additional linear inequality constraint} that \textit{cuts off} (excludes) the current solution while preserving all the original feasible solutions, and repeat  these procedures.

For CSIP, Blankenship and Falk proposed the concept of the ideal CP method \cite{blankenship1976infinitely}, but this ideal method needs to solve every subproblem without numerical errors, which is an unrealistic assumption in many CSIP problems.
In \cite{oustry2025convex}, Oustry and Cerulli studied convergence of a CP method that allows slight errors in solutions of the \textit{lower-level problems}: $\mathop{\mathrm{maximize}}_{\tau\in [0, 1]}{g_i(\bm x_r, \tau)}$ ($i=1,2,\ldots,I$), where $\bm x_r$ is the solution of the \textit{upper-level problem} $P[\widetilde{\mathcal T}_r]$ in the $r$th update.
Convergence of this CP method to the original CSIP solution is proven, including its rate, under the following assumptions.
\begin{Assu}[Strong Convexity]
In (7), the function $f$  is strongly convex, i.e., $f-\gamma\lVert\cdot\rVert^2_{2}/2$ is convex for some $\gamma>0$.
Note that it only needs to be satisfied within the feasible set.
\end{Assu}

\begin{Assu}[Linear Constraints]
In (7), every function $g_i(\cdot,\tau)$ ($i=1,2,\ldots, I$) is linear (or affine) for any $\tau\in T$,~and  every function $h_j$ ($j=1, 2, \ldots, J$) is also linear (or affine).
\end{Assu}

\begin{Assu}[Slater's Condition]
In (7), there is $\hat{\bm x} \in \mathbb R^\zeta$ s.t.~$g_i(\hat{\bm x},\tau) < 0$ ($\forall\tau \in T_i$) for all $i$ and $h_j(\hat{\bm x})= 0$ for all $j$.
\end{Assu}

To solve a CSIP problem satisfying Assumptions 1--3, the CP method in \cite{oustry2025convex} performs the following iterative procedures without approximations for the cost and constraint functions, i.e., by letting $\tilde f=f$, $\tilde g_i=g_i$, and $\tilde h_j=h_j$ in (10).

\begin{itemize}
\item[(i)]
Set a tolerance $\epsilon \ge 0$ as one input argument and~define $\widetilde T_{i,0}:=\varnothing$  ($i=1,2,\ldots,I$).
Set $r=0$ and proceed to (ii).

\item[(ii)]
Solve the upper-level problem (relaxed problem) $P[\widetilde{\mathcal{T}}_{r}]$ with $\widetilde{\mathcal{T}}_{r}:=\bigtimes^I_{i=1}\widetilde T_{i,r}$ by an appropriate  algorithm, and obtain its \textit{unique} solution $\bm x_r$.
Proceed to (iii).

\item[(iii)]
For each $i$, \textit{approximately} solve the lower-level problem
\mathindent=7mm
\begin{equation}
\mathop{\mathrm{maximize}}_{\tau\in [0, 1]}\,{g_i(\bm x_r, \tau)}\mbox{,}
\end{equation}
and obtain\footnote{If there are multiple solutions, set $\tau_{i,k}$ to the value closest to~$0$.} its solution $\tau_{i, r}$.
Proceed to (iv).

\item[(iv)]
If $g_i(\bm x_r, \tau_{i,r}) \le \epsilon$ for all $i$, terminate the algorithm and return $\bm{x}_r$ as the solution of $P[\mathcal T]$.
Otherwise, define\footnote{In \cite{oustry2025convex}, the update of $\widetilde{T}_{i, r+1}$ is written in a more general form.}
\mathindent=5mm
\begin{equation}
\widetilde{T}_{i,r+1}:=\left\{\begin{aligned}
&\widetilde{T}_{i, r}\cup\{\tau_{i,r}\}&&\mbox{if } g_i(\bm x_r, \tau_{i,r}) > \epsilon  \mbox{,}\\
&\widetilde{T}_{i, r}&&\mbox{otherwise,}
\end{aligned}\right.
\end{equation}
for each $i$, increment $r$ by 1, and proceed to (ii).
\end{itemize}

From the definition of $\widetilde{T}_{i,r+1}$ in Step (iv), $\widetilde{T}_{i, r} \subset \widetilde{T}_{i,r+1}$ holds for all $r$ in the CP method, i.e., sampling points $\tau\in T$ are added one by one to the discrete set $\widetilde{T}_i$ in (10) during the iterative computations, differently from the LRSQP method.
Each sampling point $\tau$ in $\widetilde{T}_i$ is the solution of the lower-level problem in (11), where $\tau_{i,r}$ in Step (iii) does not necessarily have to be the exact maximizer of $g_i(\bm x_r, \cdot)$.
For convergence, it is sufficient that $\tau_{i,r}\in T$ satisfies the following condition
\mathindent=7mm
\begin{equation}
\phi_i(\bm x_r)-\delta|\phi_i(\bm x_r)|
\le
g_i(\bm x_r, \tau_{i,r})
\le
\phi_i(\bm x_r)
\end{equation}
for some $\delta \in [0,1)$, where $\phi_i(\bm x) := \mathop{\mathrm{max}}_{\tau\in [0,1]} g_i(\bm x, \tau)$.
If $\tau_{i, r}$ is the exact maximizer, the condition in (13) holds with $\delta=0$.
Oustry and Cerulli proved the following theorem \cite{oustry2025convex}.

\begin{Theo}[Convergence of the Cost $f(\bm{x})$]
Suppose that Assumptions 1--3 and the condition in (13) hold.
Let $\bm x^*$~be the unique optimal solution of the CSIP problem in (7).
The solution $\bm x_r$ of $P[\widetilde{T}_r]$ in Step (ii) of the CP method satisfies 
\mathindent=7mm
\[
f(\bm x^*) - f(\bm x_r) \leq\frac{\eta}{\gamma(1-\delta)^2 (r+2)}\mbox{,}
\]
i.e., the convergence rate is $f(\bm x^*) - f(\bm x_r)=\mathcal{O}(1/r)$, where $\eta\geq0$ is some constant (see \cite[Theorem 5]{oustry2025convex} for details).
\end{Theo}

In the next section, first we show that the nonnegative spline smoothing in \eqref{cSIP} satisfies Assumptions 1--3.
Then, we provide the proposed algorithm based on the CP method.

\subsection{Proposed Nonnegative Spline Smoothing Algorithm}
In (9), as mentioned in Section 2.2, the cost function $f(\bm b)= \lVert \bm{y}-\bm{A}\bm{b}\rVert_{2}^{2}+\lambda\bm{b}^{\mathrm{T}}\bm{Q}\bm{b}$ is~strongly convex on the linear equality constraint $\bm {H}\bm b=\bm{0}$ from Appendix~A, and thus Assumption~1 holds.
The inequality constraint function $g_i(\bm b,\tau)={-}\bm g_\tau^{\mathrm T}\bm b_{\langle i \rangle}$ is linear with respect to $\bm{b}$, and thus Assumption 2 holds.
In addition, by letting $\hat{\bm b}=\bm 1$, we can obtain a constant function $s(x)=1$ for all $x\in [\xi_0, \xi_m]$ that satisfies $g_i(\hat{\bm b}, \tau)=-1<0$ and $h_j(\hat{\bm b})=\bm{h}_j^{\mathrm T}\hat{\bm b}=0$ obviously, i.e., Assumption 3~holds.

Since Assumptions 1--3 hold, we apply the CP method in Section 4.1 to the 1D nonnegative spline smoothing in (9).
Before providing the proposed algorithm, for each $r\in\mathbb N$, we define a sparse matrix $\bm G_r\in\mathbb R^{\sum_{i=1}^{m}{|\widetilde{T}_{i, r}|}\times(dm+1)}$ that satisfies
\mathindent=2.7mm
\begin{align}
\bm G_{r}\bm b
= \mathrm{vec}(\bm g_{\tau}^\mathrm{T} \bm b_{\langle i \rangle})\;\mbox{s.t.}\; \tau\in \widetilde{T}_{i,r}\ \mbox{and}\;i=1,2,\ldots, m\mbox{.}
\end{align}
By using this notation, we propose the following algorithm.

\begin{itemize}
\item[(i)]
Input nonnegative data points $\{(x_i,y_i)\}_{i=0}^{n}$ s.t.~$y_i\geq0$.
Set a tolerance $\epsilon\geq0$ as one input argument and define $\widetilde{T}_{i,0}:=\varnothing$ ($i=1,2,\ldots, m$). Set $r=0$ and proceed~to~(ii).

\item[(ii)]
By using a QP algorithm, solve the upper-level problem
\mathindent=2.8mm
\begin{equation}
\mathop{\mathrm{minimize}}_{\bm b \in \mathbb{R}^{dm+1}}\,
f(\bm b)\quad \mbox{subject to}\; \bm G_{r}\bm b \ge \bm 0\mbox{, } \bm H \bm b = \bm 0\mbox{,}\label{rbSIP}
\end{equation}
and obtain its unique solution $\bm b_r$.
Proceed to~(iii).

\item[(iii)]
For each $i$, solve the lower-level problem that is equivalent to the minimization of the polynomial piece $p_i$ in (1) with respect to $\tau$, and obtain its solution $\tau_{i, r}\in[0,1]$ through a \textit{polynomial equation}.
Proceed to (iv).

\item[(iv)]
If the minimum value of $p_i$ at $\tau_{i,r}$ given in (iii) is equal to or greater than $-\epsilon$ for all $i$, terminate the algorithm and~return $\bm{b}_r$ as the coefficients of the optimal nonnegative spline function in (4).
Otherwise, define $\widetilde{T}_{i,r+1}$ as in (12) for each $i$, increment $r$ by 1, and proceed to (ii).
\end{itemize}
We have the following theorem (see Appendix C for proof).

\begin{Theo}[Convergence of the Coefficients $\bm{b}$]
Let $\bm b^*$ be the unique optimal solution of the CSIP problem in (9).
The coefficients $\bm b_r$ in Step (ii) of the proposed algorithm satisfy
\mathindent=7mm
\[
\lVert\bm b^* - \bm b_r\rVert^2_2 \le \frac{2}{\gamma}\,\bigl(f(\bm b^*)- f(\bm b_r)\bigr)\mbox{,} 
\]
i.e., the convergence rate is at least $\lVert\bm b^* - \bm b_r\rVert^2_2=\mathcal{O}(1/r)$~under the assumption that $\tau_{i,r}$ satisfies the condition as in (13), where $\gamma$ is the positive smallest eigenvalue of the \textit{Riemannian Hessian} of $f$ on the null space $\mathcal N_{\bm{H}}:=\{\bm b\in \mathbb R^{dm+1}\,|\,\bm H \bm b=\bm 0\}$.
\end{Theo}

Differently from the LRSQP method in Section 3.2,~the proposed method directly minimizes the original cost function $f$ within the relaxed constraints in Step (ii).
This relaxed problem in (15) is a QP problem, and it is efficiently solved by an interior-point method \cite{Vanderbei-Carpenter-93,Altman-Gondzio-99}.
Hence, compared with the MATLAB LRSQP algorithm that minimizes the approximate cost functions $\tilde{f}_{r}$ by an active-set method, the proposed method is expected to solve the problem in (9) more quickly.

In Step (iii), we have to find the minimizer $\tau_{i,r}\in[0,1]$ of each polynomial piece $p_{i}$.
Since our target is the 1D~case, $\tau_{i,r}$ can be obtained with (almost) no errors by solving a univariate polynomial equation  as follows.
From (1), we~have 
\mathindent=1mm
\begin{align}
\frac{\mathrm{d}p_i}{\mathrm{d}\tau}&= \sum_{k=0}^{d-1} 
\frac{d!(b_{d(i-1)+k+1}-b_{d(i-1)+k})}{(d-k-1)!k!}\,
(1-\tau)^{d-k-1} \tau^k\nonumber \\
&=\sum_{l=0}^{d-1}\sum_{k=0}^{l}\frac{ d!(-1)^{l-k}(b_{d(i-1)+k+1}-b_{d(i-1)+k})}{(d-l-1)!(l-k)!k!}\,\tau^l
\end{align}
and find $\tau^*_i\in[0,1]$ at which the value of (16) takes $0$, where the coefficients $b_{d(i-1)+k}$ are fixed to the current solution~$\bm{b}_r$.
The minimizer $\tau_{i,r}$ of $p_i$ is either $0$, $1$, or $\tau^*_i$, and there are at most $d-1$ values\footnote{When all solutions of the polynomial equation ${\mathrm{d}p_i}/{\mathrm{d}\tau}=0$ do not belong to the interval $[0,1]$, $\tau_i^*$ no longer exists.
When $p_i$ is a polynomial of degree $0$, i.e., a constant function, we set $\tau_{i,k}=0$.} of $\tau^*_i$ since they are some of the solutions of the polynomial equation ${\mathrm{d}p_i}/{\mathrm{d}\tau}=0$ of degree $d-1$.

When $3\leq d \leq 5$, the degree of the polynomial ${\mathrm{d}p_i}/{\mathrm{d}\tau}$ is usually between $2$ and $4$.
Hence, we can \textit{exactly} obtain $\tau_{i}^*$ from the closed-form solutions of the polynomial equation, and the minimizer $\tau_{i,k}$ is given by comparing the values of~$p_i$ at $\tau=0, 1, \tau^*_i$.
In this case, the condition as in (13) holds with $\delta=0$.
When $d\ge6$, the degree of ${\mathrm{d}p_i}/{\mathrm{d}\tau}$ is equal to or greater than $5$, and there are no closed-form solutions of~the polynomial equation.
Instead, we utilize numerical solutions that can be calculated by a basic function equipped in major programming languages.
In this paper, we use the MATLAB \texttt{roots} function, where we input the coefficients in the form of (16).
Even in this case, we can expect that the condition as in (13) holds with some small $\delta\in(0,1)$, and the proposed algorithm is guaranteed to converge regarding both the cost function and the coefficient vector from Theorems 2 and 3.

Actually, $\bm b_r$ in Step (ii) is also a numerical solution~obtained by an interior-point method and contains slight errors.
Since the QP problem in (15) is not ill-conditioned and yields a sufficiently accurate $\bm b_r$, the proposed method works stably in practice.
Note that the smaller numerical errors are, the better the convergence performance is, in general.
Thus,~using the closed-form solutions in Step (iii) may yield shorter computation times compared to using numerical solutions.

\section{Numerical Experiments}
We compare 1D nonnegative spline smoothing by the conventional QP method with the sufficient condition \cite{kitahara2016thesis}, the LRSQP method usable in MATLAB \cite{matlab_optim_doc}, and the proposed CP method.
In the following, nonnegative data $\{(x_i,y_i)\}_{i=0}^{n}$ s.t.~$y_i\geq0$ are given at $x_i:=i$, and we set knots of all spline functions to $\xi_i:=x_i$.
We conducted the following numerical experiments by using MATLAB R2024a 64-bit on Microsoft Surface Laptop 7 with Snapdragon X Elite and 16 GB RAM.

\subsection{Comparison with the Conventional Methods}
\begin{figure}[t]
  \centering
  \includegraphics[height=43mm,width=70mm]{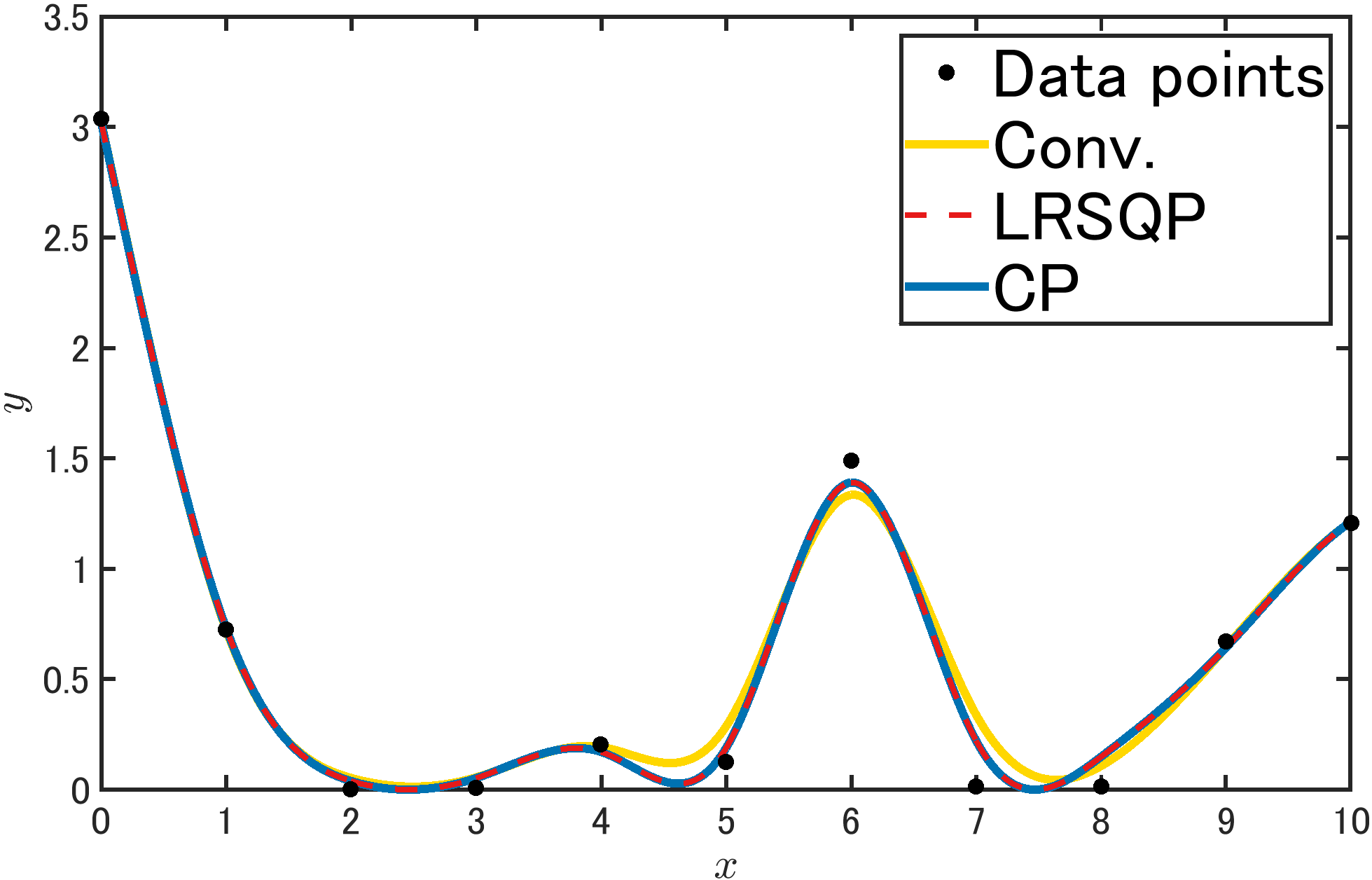}
  \caption{Nonnegative spline smoothing results for $d=3$.}
  \label{fig:d3}
\end{figure}
We compare the proposed method with the QP method and the LRSQP method in low-degree cases where $d=3,4$.
We generate nonnegative data $y_i$ as follows.
First, we generate~$\tilde y_i$ ($i=0,1,2,\ldots,n$) from the standard normal distribution of mean $0$ and variance $1$.
Then, so that nonnegative data~points very close to $0$ are repeatedly given, we define $y_i:=|\tilde y_i|/100$ if $i=5t+2$ or $i=5t+3$ for $t\in\mathbb N$, and $y_i:=|\tilde y_i|$ otherwise.
For the QP method under the sufficient condition in (5), we use the MATLAB \texttt{quadprog} function with default settings, where an interior-point algorithm once solves the problem in (6).
For the LRSQP method, we use the MATLAB \texttt{fseminf} function with default settings, where an active-set algorithm iteratively solves the approximate problems with respect to the direction vector.
For the proposed CP method, we also use MATLAB \texttt{quadprog} function with the default settings in Step (ii), i.e., the interior-point algorithm iteratively solves the relaxed problems in (15).
In Step (iii), we compute each $\tau^*_i$ by the closed-form solutions of the polynomial equation of degree $2$ or $3$.
We set the tolerance parameter used for the termination condition of the proposed method to $\epsilon = 0$.
We fix the smoothing parameter to $\lambda = 1/250$ so that results of the spline smoothing in (3) violate the overall nonnegativity.

\begin{figure}[t]
  \centering
  \includegraphics[height=43mm,width=70mm]{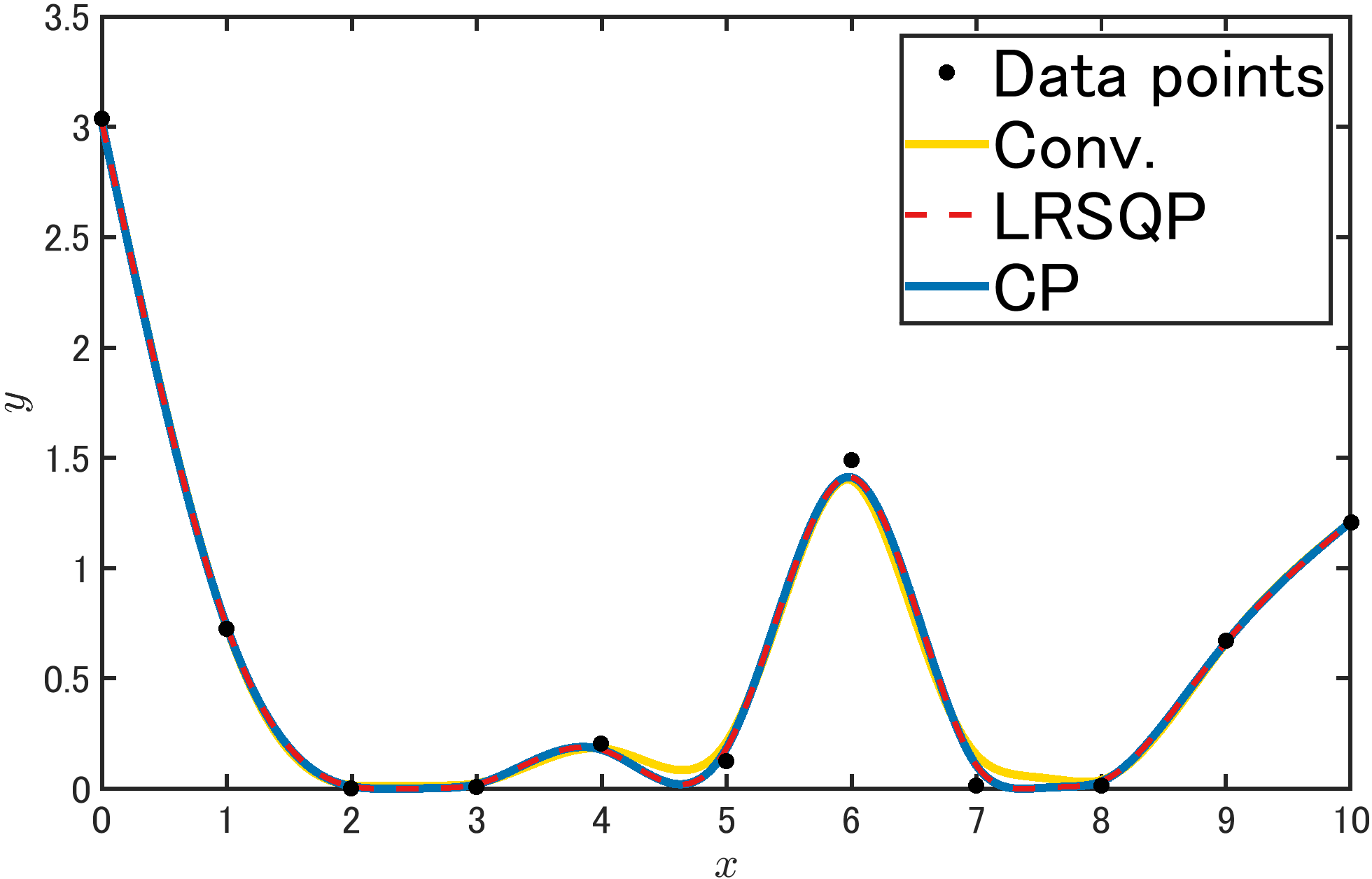}
  \caption{Nonnegative spline smoothing results for $d=4$.}
  \label{fig:d4}
\end{figure}    
\begin{figure}[t]
  \centering
  \includegraphics[height=43mm,width=70mm]{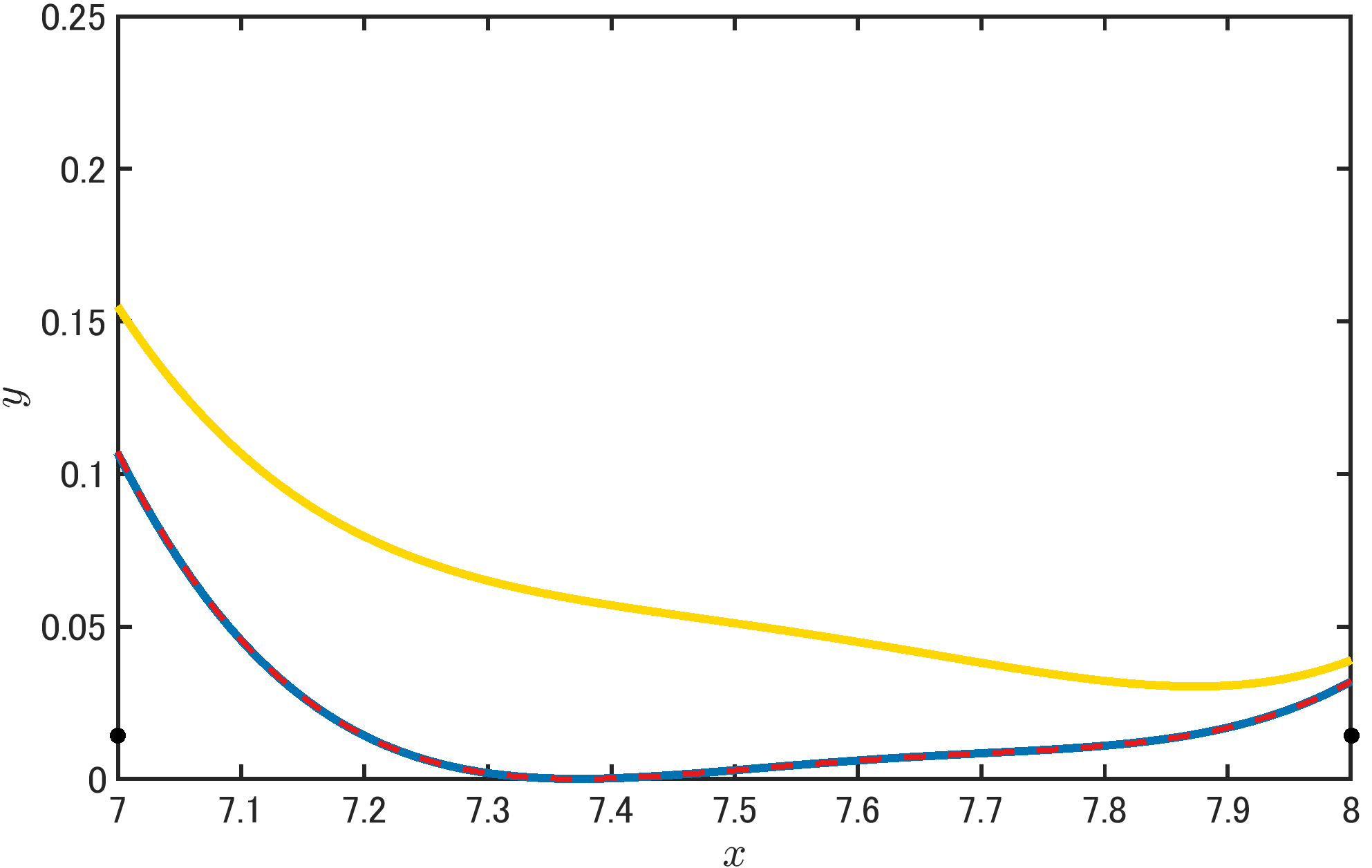}
  \caption{Magnified view of the interval $[7,8]$ in Fig.~\ref{fig:d4}.}
  \label{fig:zoom1} 
\end{figure}

Figures \ref{fig:d3} and \ref{fig:d4} show the results for the cases of $n=10$ with $d=3$ and $d=4$, respectively.
In both figures, black dots indicate the given data points $(x_i, y_i)$, and yellow, red, blue lines represent the nonnegative spline functions obtained by the QP method, the LRSQP method, and the proposed CP method, respectively.
From Fig.~1, we can observe that the results of the LRSQP method and the CP method are visually the same since both methods solve the same CSIP problem in (9).
The result of the QP method is similar in many areas, but is slightly shifted upward around the intervals $[4,5]$ and $[7,8]$, which means that the solution is degraded because of the feasible set restricted by the sufficient condition.
Figure~2 displays similar characteristics to Fig.~1, but by changing the degree $d$ from $3$ to $4$, the difference between the QP method and the two CSIP methods becomes smaller.
We can see~that increasing the degree $d$ improves the flexibility, resulting in the nonnegative spline functions that lie closer to each data.
This effect tends to be particularly notable in areas where the function values approach zero such as the interval $[7,8]$.

Figure \ref{fig:zoom1} shows a magnified view of the interval $[7,8]$ in Fig.~2.
As mentioned above, the result of the QP method lies on the upper side of the results of the two CSIP methods, i.e., the sufficient condition for the nonnegativity is preventing the spline function from approaching zero.
On the other hand, the LRSQP method and the proposed CP method obtain the nonnegative smooth curves that naturally pass very close to zero.
The minimum values of the spline functions in~$[7,8]$ are $3.028\times10^{-2}$ for the QP method, $1.505\times10^{-4}$ for the LRSQP method, and $7.752\times10^{-10}$ for the CP method.
From these values, we see that both CSIP methods can construct nonnegative spline functions that are closer to zero than the QP method, and the proposed CP method obtains the most accurate spline function whose minimum value is almost 0.

\begin{table}[t]
\centering
\caption{Comparison of the cost function values at the solutions.}
\resizebox{\linewidth}{!}{
\begin{tabular}{c ccc ccc}
\toprule
 & \multicolumn{3}{c}{$d=3$} & \multicolumn{3}{c}{$d=4$} \\
\cmidrule(lr){2-4} \cmidrule(lr){5-7}
$n$ & QP & LRSQP & CP & QP & LRSQP & CP \\
\midrule
5   & 0.0210730 & 0.0190703 & \bf{0.0190702} & 0.0174357 & 0.0168465 & \bf{0.0168453} \\
10  & 0.2416261 & \bf{0.1752847} & \bf{0.1752847} & 0.1494654 & 0.1339960 & \bf{0.1339937} \\
50  & 1.2113624 & \bf{0.8621099} & \bf{0.8621099} & 0.7282246 & 0.6468150 & \bf{0.6468088} \\
100 & 2.4235341 & \bf{1.7206420} & \bf{1.7206420} & 1.4516748 & 1.2878329 & \bf{1.2878277} \\
\bottomrule
\end{tabular}}
\label{tab:cost}
\end{table}
\begin{table}[t]
\centering
\caption{Comparison of the computation times (seconds).}
\resizebox{0.9\linewidth}{!}{
\begin{tabular}{c ccc ccc}
\toprule
 & \multicolumn{3}{c}{$d=3$} & \multicolumn{3}{c}{$d=4$} \\
\cmidrule(lr){2-4} \cmidrule(lr){5-7}
$n$ & QP & LRSQP & CP & QP & LRSQP & CP \\
\midrule
5   & 0.0077 & 0.1043 & 0.0225 & 0.0079 & 0.1375 & 0.0346 \\
10  & 0.0105 & 0.2718 & 0.0277 & 0.0126 & 0.4345 & 0.0433 \\
50  & 0.0206 & 4.7286 & 0.2453 & 0.0217 & 9.9481 & 0.4514 \\
100 & 0.0293 & 20.131 & 0.8601 & 0.0340 & 45.586 & 1.0463 \\
\bottomrule
\end{tabular}}
\label{tab:time}
\end{table}

Table \ref{tab:cost} shows the cost function value $f(\bm{b})$ at the solution obtained by each method for different $n=5,10,50,100$, where the data points for $n=10$ are the same as in Figs.~1--3.
From Table~1, the two CSIP methods reduced the cost~function values by up to about $30\%$ for $d=3$ and $10\%$ for $d=4$ from the QP method with the sufficient condition.
Compared to the LRSQP method, the cost function values by the~proposed CP method are almost the same for $d=3$, and slightly better for $d=4$.
Note that the cost function values for $d=4$ are smaller than those for $d=3$ in all the methods.
This is because $\mathcal S^{2}_{d}(\sqcup_m)\subset\mathcal S^{2}_{d+1}(\sqcup_m)$ holds for any partition $\sqcup_m$~and $d\geq3$, and increasing the degree $d$ results in an equivalent or superior solution along with a longer computation time.

Table \ref{tab:time} shows the computation times of each method.
From Table \ref{tab:time}, the QP method under the sufficient condition requires the shortest computation time because this method solves the problem in (6) only once, but its solution deteriorates from the original CSIP solution as shown in Figs.~1--3 and Table 1.
Comparing the two CSIP methods, the proposed CP method requires less computation time than the~LRSQP method in all the cases, which may be because the convergence performance of the LRSQP method is degraded by~the inner iterative approximations including numerical differentiation.
The difference is significant when the number of~data points is large especially for $d=4$.
For $n=100$, the computation time of the proposed method was about $4\%$ of that of the LRSQP method for $d=3$, and was about $2\%$ for $d=4$.

We found that the minimum values of the nonnegative spline functions computed by the proposed CP method were very slightly greater than zero.
If the result of the standard spline smoothing in (3) does not satisfy the overall nonnegativity, then the exact optimal solution of the nonnegative spline smoothing in (9) lies on the boundary of the feasible set, i.e., its minimum value is exactly $0$.
Thus, the proposed method with the interior-point algorithm in Step (ii) can obtain a nonnegative spline function that lies \textit{slightly inside} the exact solution.
This property is desirable in practice because obtaining a spline function that lies \textit{slightly outside} the exact solution requires some post-processing to cancel out slight negative function values.
From the above, we confirmed that the proposed CP method can quickly compute the solution of the nonnegative spline smoothing with very high accuracy.

\subsection{Effects of the Degree of Spline Functions}
We next investigate the behavior of the proposed CP method for higher degrees $d=5,6,\ldots,10$, where candidates of~$\tau_i^*$ in Step (iii) were computed from numerical solutions of the polynomial equation by using the MATLAB \texttt{roots} function.
The number of data points was fixed to $n = 10$, i.e., the data points are the same as in Figs.~1--3, with the same tolerance and smoothing parameters $\epsilon=0$ and $\lambda=1/250$.
To compare with the results for the low-degree cases in Tables~1 and 2, we applied the proposed method with the MATLAB \texttt{roots} function not only for $d=5,6,\ldots,10$ but also for $d=3,4$.
\begin{figure}[t]
  \centering
  \includegraphics[height=43mm,width=70mm]{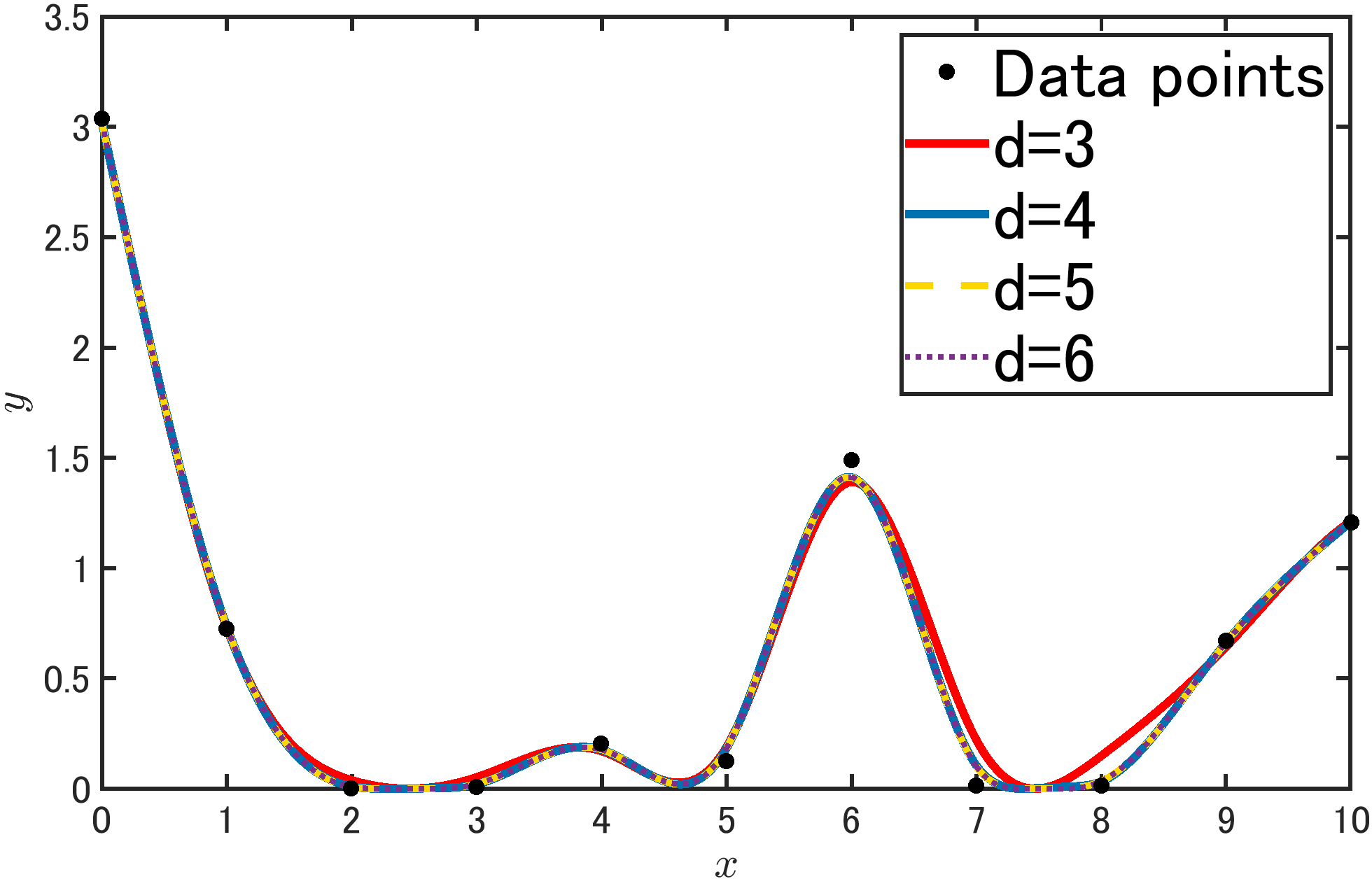}
  \caption{Nonnegative spline smoothing results for $d=3,4,5,6$.}
  \label{fig:deg}
\end{figure}

\begin{table}[t]
\centering
\caption{Convergence performance of the CP method for each degree.}
\resizebox{\linewidth}{!}{\begin{tabular}{c cccc}
\toprule
$d$ & 3 & 4 & 5 & 6 \\
\midrule
Cost Function Value
& 0.1752847 & 0.1339937 & 0.1336444 & 0.1335042 \\
Computation Time $[\mathrm{s}]$
& 0.0419 & 0.0553 & 0.0607 & 0.0927 \\
\midrule
$d$ & 7 & 8 & 9 & 10 \\
\midrule
Cost Function Value
& 0.1334590 & 0.1334577 & 0.1334472 & 0.1334441 \\
Computation Time $[\mathrm{s}]$
& 0.0961 & 0.1068 & 0.1103 & 0.1778 \\
\bottomrule
\end{tabular}}
\label{tab:degree}
\end{table}

Figure \ref{fig:deg} shows the results of the proposed CP method with the MATLAB \texttt{roots} function, where red, blue, yellow, and purple lines respectively correspond to $d=3$, $4$, $5$, and $6$.
Table \ref{tab:degree} summarizes the cost function values at the solutions and the computation times for $d=3,4,\ldots,10$.
From Fig.~4, the results for $d=3$ and $d=4$ are visually distinct as in the comparison of Figs.~1 and 2, and their cost function values in Table~3 are the same as in Table~1.
On the other hand, the results for $d=4,5,6$ are visually identical in Fig.~4.
In fact, the cost function values in Table~3 change only slightly for $d \ge 4$, although the flexibility of spline functions increases.

Comparing the computation times in Table~2 with those in Table~3 for $d=3,4$, we see that the using the MATLAB \texttt{roots} function in Step (iii) instead of the closed-form solutions increased the computation time by about $85\%$ for~$d=3$ and $60\%$ for $d=4$.
Increasing the degree $d$ further~extended the computation time in Table~3.
These results indicate that, even for higher-degree cases, the proposed method with the MATLAB \texttt{roots} function correctly obtains the solution, but setting the degree to $d\geq5$ provides little benefit in practice.

From the results in Sections 5.1 and 5.2, we recommend to use the proposed CP method for $d=4$, together~with~an interior-point~algorithm in Step (ii) and the closed-form solutions in Step~(iii), which is the best in terms of the balance between speed and accuracy.
If a shorter computation time is needed, one possible approach is to set the tolerance parameter $\epsilon$ to a small positive number.
In this case, the computation~time will be shorter, but a spline function having slight negative values will be obtained as the solution.
By~adding a positive constant function (i.e., a positive polynomial of degree $0$)~to cancel out the negative function values as post-processing, we can construct a nonnegative spline function close to the exact solution in a shorter computation time.

\section{Conclusion}
In this paper, we newly formulated the 1D nonnegative spline smoothing as a CSIP problem, and proposed a CP method~to quickly and accurately solve this CSIP problem.
We~proved that the proposed CP method converges to the optimal spline function by using the closed-form solutions or numerical solutions of polynomial equations.
The numerical experiments showed that, compared to the conventional QP method under the sufficient condition for the nonnegativity, the proposed method constructs better nonnegative spline functions that achieve lower cost function values.
In addition,  compared to the LRSQP method, the proposed method with~an interior-point algorithm can greatly reduce computation time while yielding a slightly more accurate solution.
From a practical view, in the proposed method, setting the degree to $d=4$ and using the closed-form solutions provide a better balance between computation time and numerical accuracy~than~setting the degree to $d\geq5$ and using the numerical solutions.

\appendix
\section{Strong Convexity of the Cost Function}
For QP in (3), its cost function $f(\bm b)= \lVert \bm{y}-\bm{A}\bm{b}\rVert_{2}^{2}+\lambda\bm{b}^{\mathrm{T}}\bm{Q}\bm{b}=\bm{b}^{\mathrm{T}}(\bm{A}^{\mathrm{T}}\bm{A}+\lambda\bm{Q})\bm{b}-2\bm{y}^{\mathrm{T}}\bm{A}\bm{b}+\bm{y}^{\mathrm{T}}\bm{y}$ is twice differentiable.
Thus, the cost function $f$ is $\gamma$\textit{-strongly convex}~on the constraint set $\mathcal N_{\bm{H}}=\{\bm b\in \mathbb R^{dm+1}\,|\,\bm H \bm b=\bm 0\}$ if and only if $\nabla^2 f(\bm b)- \gamma \bm{E}\in\mathbb R^{(dm+1)\times(dm+1)}$ with some $\gamma>0$ is \textit{positive semidefinite} for all $\bm{b}\in\mathcal N_{\bm{H}}$, where $\nabla^2 f(\bm b)=2(\bm A^{\mathrm T}\bm A+\lambda\bm Q)$ is the \textit{Hessian}~of $f$ and $\bm{E}$ is the identity matrix.
Thus, we need to show~that
\mathindent = 0.5mm
\begin{equation}
\quad2\bm{b}^{\mathrm{T}}\bm A^{\mathrm T}\bm A\bm{b}+2\lambda\bm{b}^{\mathrm{T}}\bm Q\bm{b}\geq\gamma\lVert\bm{b}\rVert^{2}_{2} \quad\mbox{for all}\;\bm{b}\in\mathcal N_{\bm{H}}\label{a1}
\end{equation}
holds for some $\gamma>0$, to prove the strong convexity of $f$.

First of all, for each vector $\bm{b}\in\mathcal N_{\bm{H}}$, there is exactly one corresponding spline function $s\in\mathcal S^{2}_{d}(\sqcup_m)$. 
For the first~term in (\ref{a1}), we have $\bm{b}^{\mathrm{T}}\bm A^{\mathrm T}\bm A\bm{b}=\lVert \bm{A}\bm{b}\rVert_{2}^{2}=\sum_{i=0}^{n}|s(x_i)|^2\geq0$,  where the equality holds only when $s(x_i) = 0$ holds for all $i$.
For the second term, we have $\bm b^{\mathrm T}\bm Q\bm b=\int_{\xi_0}^{\xi_m}|s''(x)|^2\,\mathrm{d}x\geq0$, where the equality holds when the spline function is a single polynomial of degree $1$ at most on the entire interval $[\xi_0,\xi_m]$, i.e., $s(x)=c_1 x+c_0$ holds for all $x\in[\xi_0,\xi_m]$ with $c_0,c_1\in\mathbb R$.
Thus, a spline function $s$ satisfying $\bm{b}^{\mathrm{T}}\bm A^{\mathrm T}\bm A\bm{b}+\lambda\bm{b}^{\mathrm{T}}\bm Q\bm{b} = 0$ with $\lambda>0$ is only the zero function $s\equiv0$, and its coefficient vector is $\bm{b}=\bm{0}$.
Consequently, $\bm A^{\mathrm T}\bm A+\lambda\bm Q$ is \textit{positive definite} on $\mathcal N_{\bm{H}}$, i.e., $\bm{b}^{\mathrm{T}}\bm A^{\mathrm T}\bm A\bm{b}+\lambda\bm{b}^{\mathrm{T}}\bm Q\bm{b}>0$ for all $\bm{b}\in\mathcal N_{\bm{H}}\setminus\{\bm{0}\}$.

The maximum value of $\gamma$ in (\ref{a1}) is given by solving the following \textit{Rayleigh quotient minimization on} $\mathcal N_{\bm{H}}$ as
\mathindent = 7mm
\begin{equation}
\frac{\gamma}{2}=\min_{\bm{b}\in\mathcal N_{\bm{H}}\setminus\{\bm{0}\}}\frac{\bm{b}^{\mathrm{T}}(\bm A^{\mathrm T}\bm A+\lambda\bm Q)\bm{b}}{\bm{b}^{\mathrm{T}}\bm{b}}>0\mbox{.}\label{a2}
\end{equation}
By defining a matrix $\bm{V}\in\mathbb R^{(dm+1)\times(dm-2m+3)}$ whose column vectors $\bm{v}_{l}$ s.t.~$\lVert\bm{v}_{l}\rVert_{2}=1$ ($l=1,2,\ldots,dm-2m+3$) form an \textit{orthogonal basis} of the null space $\mathcal N_{\bm{H}}$, each $\bm{b}\in\mathcal N_{\bm{H}}\setminus\{\bm{0}\}$ is expressed as $\bm{b}=\bm{V}\bm{z}$ with some vector $\bm{z}\in\mathbb R^{dm-2m+3}\setminus\{\bm{0}\}$.
From this expression, we have $\bm{b}^{\mathrm{T}}\bm{b}=\bm{z}^{\mathrm{T}}\bm{V}^{\mathrm{T}}\bm{V}\bm{z}=\bm{z}^{\mathrm{T}}\bm{z}$ and
\begin{equation}
\frac{\bm{b}^{\mathrm{T}}(\bm A^{\mathrm T}\bm A+\lambda\bm Q)\bm{b}}{\bm{b}^{\mathrm{T}}\bm{b}}
=
\frac{\bm{z}^{\mathrm{T}}\bm{V}^{\mathrm{T}}(\bm A^{\mathrm T}\bm A+\lambda\bm Q)\bm{V}\bm{z}}{\bm{z}^{\mathrm{T}}\bm{z}}\mbox{,}
\label{a3}
\end{equation}
and the minimum value of the right-hand side of (\ref{a3}) is the smallest eigenvalue of  the matrix $\bm{V}^{\mathrm{T}}(\bm A^{\mathrm T}\bm A+\lambda\bm Q)\bm{V}$.
This is a positive value from (\ref{a2}), and twice this value is $\gamma$ satisfying the condition in (\ref{a1}).
Thus, $f$ is strongly convex on $\mathcal N_{\bm{H}}$.

\section{KKT Conditions for the Problem in (9)}
For simplicity of the notation of the KKT conditions, we~use a matrix $\bm G_{r+1}\in\mathbb R^{\sum_{i=1}^{m}{|\widetilde{T}_{i,r+1}|}\times(dm+1)}$ defined in the same~manner as in (14), where assume that every local maximizer in $\widetilde{T}_{i,r+1}$ has been exactly detected without numerical errors for each $g_i(\bm {b}_r, \cdot)$ in Step (iii).
If the current solution (coefficient vector) $\bm{b}_r\in\mathbb R^{dm+1}$ satisfies the following KKT conditions
\begin{equation}
\left\{
\begin{aligned}
&\nabla f(\bm{b}_r) -\bm G_{r+1}^{\mathrm T} \bm \mu + \bm H^{\mathrm T} \bm \nu= \bm 0\mbox{,}\\
&\bm G_{r+1}\bm b_r\geq\bm{0}\mbox{,}\quad\bm H \bm b_r = \bm 0\mbox{,}\\
&\bm{\mu}\geq\bm{0}\mbox{,}\quad\bm{\mu}\odot(\bm G_{r+1}\bm b_r) = \bm{0}\mbox{,}
\end{aligned}
\right.
\label{a4}
\end{equation}
with Lagrange multipliers $\bm{\mu}\in\mathbb R^{\sum_{i=1}^{m}|\widetilde{T}_{i,r+1}|}$ and $\bm{\nu}\in\mathbb R^{2m-2}$, then $\bm {b}_r$ is the optimal solution of the CSIP problem in (9), where $\odot$ is the Hadamard (component-wise) product.
In the LRSQP method, the gradient $\nabla f(\bm{b})=2\bm A^{\mathrm T}(\bm A \bm b - \bm y)+ 2 \lambda \bm Q \bm b$ of $f$ in (\ref{a4}) is approximated by numerical differentiation.

The reason why the conditions for the optimal solution of the CSIP problem including the infinite constraints can~be expressed in the finite terms lies in the last equation in (\ref{a4}).
If we consider the KKT conditions for a problem including a nearly infinite number of inequality constraints, the number of Lagrange multipliers contained in $\bm{\mu}$ is also nearly infinite, i.e., $\bm{\mu}$ behaves more like a function than a vector.
Here, by utilizing the last equation of the KKT conditions, we see that if $g_i(\bm{b}_r, \tau)<0$, we can exclude the corresponding Lagrange multiplier from the conditions since its value is $0$.
Therefore, as in (\ref{a4}), it is sufficient to check the KKT conditions only for the discrete set $\widetilde{T}_{i,r+1}$ consisting of the local maximizers (or consisting of the global maximizer as in the CP method).

\section{Proof of Theorem 3 with General Cases}
From the $\gamma$-strong convexity of $f(\bm b)= \lVert \bm{y}-\bm{A}\bm{b}\rVert_{2}^{2}+\lambda\bm{b}^{\mathrm{T}}\bm{Q}\bm{b}$ on the null space $\mathcal N_{\bm{H}}=\{\bm b\in \mathbb R^{dm+1}\,|\,\bm H \bm b=\bm 0\}$ of $\bm H$, we~have
\mathindent=4mm
\begin{equation}
f(\bar{\bm b})\ge f(\bm b)+\langle\mathrm{grad}\,f(\bm b),{\bar{\bm b}-\bm b}\rangle+\frac{\gamma}{2}\lVert\bar{\bm b}-\bm b\rVert_2^2
\label{a5}
\end{equation}
for any $\bm b,\bar{\bm b}\in\mathcal N_{\bm{H}}$, where $\langle\cdot , \cdot \rangle$ is the (standard) inner product, and the standard gradient $\nabla f(\bm{b})=2\bm A^{\mathrm T}(\bm A \bm b - \bm y)+ 2 \lambda \bm Q \bm b$ is replaced with the \textit{Riemannian gradient} $\mathrm {grad}\, f(\bm b)\in \mathbb R^{dm+1}$ on $\mathcal N_{\bm{H}}$ \cite{Boumal-23} by using $\langle\nabla f(\bm b),{\bar{\bm b}-\bm b}\rangle=\langle\mathrm{grad}\,f(\bm b),{\bar{\bm b}-\bm b}\rangle$.
The Riemannian gradient in (\ref{a5}) is specifically 
given by 
\mathindent=1.9mm
\[
\mathrm {grad}\, f(\bm b) = \bm P_{\mathcal N_{\bm{H}}} \nabla f(\bm b)=\bigl(\bm E-\bm H^\mathrm{T} (\bm H\bm H^\mathrm{T})^{-1}\bm H\bigr)\nabla f(\bm b)\mbox{,}
\]
where the matrix $\bm P_{\mathcal N_{\bm{H}}}$ is the orthogonal projection onto $\mathcal N_{\bm{H}}$.
By substituting $\bm{b}_{r}$ into $\bm{b}$ and $\bm{b}^{*}$ into $\bar{\bm{b}}$ in (\ref{a5}), we have
\mathindent=0.6mm
\begin{equation}
\lVert{\bm b^*}-\bm b_r\rVert_2^2\leq
\frac{2}{\gamma}\,\bigr(f({\bm b^*})-
f(\bm b_r)
-
\langle\mathrm{grad}\, f(\bm b_r),{{\bm b^*}-\bm b_r}\rangle\bigr)\mbox{.}\label{a6}
\end{equation}
We show that $\langle\mathrm{grad}\, f(\bm b_r),{{\bm b^*}-\bm b_r}\rangle\geq0$ holds in (\ref{a6}), to prove the equation in Theorem 3.

In the following, we show $\langle\mathrm{grad}\, f(\bm b_r),{{\bm b^*}-\bm b_r}\rangle\geq0$~for general constraint functions $g_{i}$ ($i=1,2,\ldots,I$), where each $g_i(\cdot,\tau)$ is differentiable and convex with respect to $\bm b$ for any $\tau\in T$.
The vector $\bm b_r$ in Step (ii) is the solution of the relaxed problem and satisfies the KKT conditions.
Hence, we have
\mathindent=6.2mm
\begin{equation}
\nabla f(\bm b_r)
+\sum_{i=1}^{I}\sum_{\kappa=1}^{K_i}\mu_{i,\kappa}\nabla g(\bm b_r,\tau_{i,\kappa})+\bm H^\mathrm{T}\bm \nu
=\bm 0\label{a7}
\end{equation}
with Lagrange multipliers $\mu_{i,\kappa}\geq 0$ and $\bm{\nu}\in\mathbb R^{2m-2}$,~where we suppose $\widetilde{T}_{i,r}=\{\tau_{i,1},\tau_{i,2},\ldots,\tau_{i,K_i}\}$ s.t.~$K_i\leq r$ and $\nabla g$ is the (standard) gradient of $g(\cdot,\tau)$ with respect to $\bm b$ for a fixed $\tau\in T$.
For the problem in (15), the second term in (\ref{a7}) is expressed as $-\bm G_r^\mathrm{\mathrm{T}} \bm \mu$.
By applying the orthogonal projection $\bm P_{\mathcal N_{\bm{H}}}$ to (\ref{a7}) from the left side, we have
\mathindent=7mm
\begin{equation}
\mathrm{grad}\, f(\bm b_r)=
-\sum_{i=1}^{I}\sum_{\kappa=1}^{K_i}\mu_{i,\kappa}\,\mathrm{grad}\, g(\bm b_r,\tau_{i,\kappa})\mbox{.}
\label{a8}
\end{equation}
From (\ref{a8}), we only need to prove $\langle\mathrm{grad}\, f(\bm b_r),\bm b^*-\bm b_r\rangle
=
-\sum_{i=1}^{I}\sum_{\kappa=1}^{K_i}\mu_{i,\kappa}
\langle\mathrm{grad}\, g_i(\bm b_r,\tau_{i,\kappa}),\bm b^*-\bm b_r\rangle\geq0$.
For this, we show $-\langle\mathrm{grad}\, g_i(\bm b_r,\tau_{i,\kappa}), \bm b^*-\bm b_r\rangle\geq g_i(\bm b_r,\tau_{i,\kappa})$.
Since $g_i(\cdot,\tau_{i,\kappa})$ is convex with respect to $\bm{b}$ for all $i$ and $\kappa$, we have
\mathindent=3mm
\[
g_i(\bm b^*,\tau_{i,\kappa})\geq g_i(\bm b_r,\tau_{i,\kappa})+\langle\mathrm{grad}\, g_i(\bm b_r,\tau_{i,\kappa}), \bm b^*-\bm b_r\rangle\mbox{.}
\]
Note that, for the problem in (15), this inequality becomes~the equality since $g_{i}(\cdot,\tau_{i,\kappa})$ is linear.
Then, we have
\mathindent=2.5mm
\begin{align*}
-\langle\mathrm{grad}\, g_i(\bm b_r,\tau_{i,\kappa}), \bm b^*-\bm b_r\rangle
&\geq
g_i(\bm b_r,\tau_{i,\kappa})-g_i(\bm b^*,\tau_{i,\kappa})\nonumber\\
&\geq g_i(\bm b_r,\tau_{i,\kappa})
\end{align*}
because the optimal solution $\bm b^*$ of the problem in (9) satisfies $g_i(\bm b^*,\tau)\le0$ for all $\tau\in T$.
Thus, recalling $\mu_{i,\kappa}\geq0$, we~have
\mathindent=2mm
\begin{equation}
\langle\mathrm{grad}\, f(\bm b_r),\bm b^*-\bm b_r\rangle
\ge
\sum_{i=1}^{I}\sum_{\kappa=1}^{K_i}\mu_{i,\kappa} g_i(\bm b_r,\tau_{i,\kappa})=0
\label{a9}
\end{equation}
from the complementary slackness $\mu_{i,\kappa} g_i(\bm b_r,\tau_{i,\kappa})=0$ for all $i$ and $\kappa$.
From (\ref{a6}) and (\ref{a9}), we obtain the equation in Theorem 3.
Moreover, if the condition as in (13) is satisfied, then by using the equation in Theorem 2, we have
\mathindent=2.7mm
\[
\lVert\bm b^* - \bm b_r\rVert^2_2 \le \frac{2}{\gamma}\,\bigl(f(\bm b^*)- f(\bm b_r)\bigr)\leq\frac{2\eta}{\gamma^2(1-\delta)^2 (r+2)}\mbox{,}
\]
i.e., the convergence rate of $\lVert\bm b^* - \bm b_r\rVert^2_2$ is at least $\mathcal{O}(1/r)$.

Finally, the value of $\gamma>0$ is the smallest eigenvalue of $2\bm{V}^{\mathrm{T}}(\bm A^{\mathrm T}\bm A+\lambda\bm Q)\bm{V}=\bm{V}^{\mathrm{T}}\nabla f^2(\bm b)\bm{V}$ as shown in Appendix~A, while the \textit{Riemannian Hessian} $\mathrm {Hess}\, f(\bm b)\in\mathbb R^{(dm+1)\times(dm+1)}$ of $f$ on $\mathcal N_{\bm{H}}$ \cite{Boumal-23} is specifically given as a fixed matrix by
\mathindent=4.3mm
\begin{align*}
\mathrm {Hess}\, f(\bm b)& = 2\bm P_{\mathcal N_{\bm{H}}}(\bm A^{\mathrm T}\bm A+\lambda\bm Q)\bm P_{\mathcal N_{\bm{H}}}\nonumber\\
&=
\bm P_{\mathcal N_{\bm{H}}} \nabla f^2(\bm b)\bm P_{\mathcal N_{\bm{H}}}=\bm{V}\bigl(\bm{V}^{\mathrm{T}}\nabla f^2(\bm b)\bm{V}\bigr)\bm{V}^{\mathrm{T}}\mbox{,}
\end{align*}
where we use the relation $\bm P_{\mathcal N_{\bm{H}}} =\bm{V}\bm{V}^{\mathrm{T}}$.
Applying $\bm{V}$ and $\bm{V}^{\mathrm{T}}$ to both sides of $\bm{V}^{\mathrm{T}}\nabla f^2(\bm b)\bm{V}$ increases the size of the matrix and adds zero eigenvalues, but does not change the values of the positive eigenvalues.
Therefore, $\gamma$ is the positive smallest eigenvalue of the Riemannian Hessian $\mathrm {Hess}\, f(\bm b)$ on $\mathcal N_{\bm{H}}$.
\end{document}